\documentclass[12pt]{article}%
\usepackage{amsmath}
\usepackage{amsfonts}
\usepackage{amssymb}
\usepackage{graphicx}
\usepackage{amsbsy}
\usepackage{amscd}
\usepackage{acronym}%
\newtheorem{theorem}{Theorem}

\newtheorem{lemma}[theorem]{Lemma}

\newenvironment{proof}[1][Proof]{\noindent\textbf{#1.} }{\ \rule{0.5em}{0.5em}}
\begin{document}

\title{On Weighted $L^{2}$ Cohomology}
\author{John G. Miller\\Indiana U.-Purdue U. at Indianapolis}
\date{}
\maketitle

\begin{abstract}
Consider an orientable manifold with countably many complete components of
bounded dimension. Suppose that its rational homology is infinitely generated
in some degree. Then there is no choice of weight function for which the
natural map from weighted $L^{2}$ cohomology to de Rham cohomology is
surjective in that degree.

\end{abstract}

\bigskip A weighted $L^{2}$ space on a Riemannian manifold $M$ is obtained by
replacing the volume measure $dm$ by a measure of the form $\phi^{2}dm,$ where
$\phi$ is a positive function referred to as a weight function. If $M$ is
noncompact and $\phi$ is sufficiently unbounded above or unbounded away from
zero, the domains of differential operators on these spaces may differ from
those on the standard $L^{2}$ space. The use of such spaces in connection with
partial differential equations has a long history. We are concerned here with
a question in $L^{2}$ cohomology. There is a natural homomorphism from
unreduced $L^{2}$ cohomology computed on the weighted spaces of differential
forms, $H_{\phi}^{\ast},$ to de Rham cohomology $H_{dR}^{\ast}.$ For which
manifolds does there exist a weight function such that this map is an
isomorphism? The first results of this sort were apparently proved by Borel
\cite{Bor} and Zucker \cite{Zuc}. Further positive results have been obtained
by a number of authors. See in particular \cite{Bue}, \cite{Bul}, \cite{Yeg},
\cite{Mil}.

There is a Hodge Laplacian $\Delta_{\phi}$ (which is $D_{\phi}^{2}$ as defined
below). Elements of its kernel in degree $k$, $\mathcal{H}_{\phi}^{k},$ are
called $\phi$-harmonic. There are injections $\mathcal{H}_{\phi}%
^{k}\rightarrow$ $H_{\phi}^{k}$ which are isomorphisms if and only if
$\Delta_{\phi}^{k}$ has closed range. Bueler formulated a general conjecture
\cite{Bue}: let $M$ be complete, connected, oriented and with Ricci curvature
bounded below. Let $\phi$ be a fundamental solution of the scalar heat
equation. Then $\Delta_{\phi}$ has closed range and $\mathcal{H}_{\phi}^{\ast
}$ is isomorphic to $H_{dR}^{\ast}.$ Carron \cite{Car} has recently found
examples which disprove this conjecture, and that do much more. Let $S$ be a
compact orientable surface of genus $\geq2$, let $\tilde{S}$ be an infinite
cyclic covering, and let $T^{n-2}$ be a torus. Let $M=\tilde{S}\times
T^{n-2}.$ There is no $\phi$ such that the natural map $\mathcal{H}_{\phi}%
^{k}\rightarrow H_{dR}^{k}$ is surjective in any degree $k$ with $0<k<n.$
Therefore, either $\Delta_{\phi}^{k}$ does not have closed range or $H_{\phi
}^{k}\rightarrow H_{dR}^{k}$ is not surjective. The present paper goes one
step further. Let $M$ be any orientable manifold with countably many complete
components of bounded dimension and $H_{k}\left(  M;\mathbb{R}\right)  $
infinite dimensional in any degree $k$. Then there is no $\phi$ such that
$H_{\phi}^{k}\rightarrow H_{dR}^{k}$ is surjective. No bounded geometry
hypothesis is required. The proof is not related to Carron's. It is shown that
$M$ may be replaced by a union of tubular neighborhoods of submanifolds
representing a basis for $H_{k}\left(  M;\mathbb{Q}\right)  .$ In this
situation the result amounts to the fact that $\ell^{2}$ doesn't contain all
sequences of numbers. It is still possible that Bueler's conjecture holds, or
that some other choice of weight produces an isomorphism, for manifolds with
finitely generated homology. This paper makes essential use of ideas of
Dodziuk \cite[Section 3]{Dod}.

We describe the analytic framework. See \cite{Bue} for more information on
this material.

$M$ is an oriented Riemannian manifold with countably many complete components
$M_{p}$ of bounded dimensions $n_{p}.$ Let $n=\max\left\{  n_{p}\right\}  .$

$\Omega_{c}^{k}$ is the complex valued smooth $k$-forms on $M$ with compact
support. If $\left\langle u,v\right\rangle $ is the standard pointwise inner
product of $k$-forms, $\left(  u,v\right)  =\int_{M}\left\langle
u,v\right\rangle dm.$ The associated $L^{2}$ norm will be written $\left\Vert
u\right\Vert .$

$\phi$ is a positive smooth weight function.

$\Omega_{\phi,c}^{k}$ is $\Omega_{c}^{k}$ with inner product $\left(
u,v\right)  _{\phi}=\int_{M}\left\langle u,v\right\rangle \phi^{2}dm=\left(
\phi u,\phi v\right)  $ and norm $\left\Vert u\right\Vert _{\phi}.$

$d^{k}$ is the exterior derivative on $\Omega_{c}^{k},$ $\delta^{k+1}$ is its
formal adjoint with respect to $\left(  \cdot,\cdot\right)  ,$ and
$\delta_{\phi}^{k+1}$ is its formal adjoint with respect to $\left(
\cdot,\cdot\right)  _{\phi}.$

$D_{\phi}=d+\delta_{\phi}$ acting on $\Omega_{\phi,c}^{\ast}.$

The closure of an operator $T$ will be written as $\bar{T}.$ The domain of
$\overline{D}_{\phi},$ $\mathcal{D}\left(  \overline{D}_{\phi}\right)  $ is
the completion of $\Omega_{\phi,c}^{k}$ for the graph norm $\left\Vert
u\right\Vert _{D_{\phi},\phi}=\left\Vert D_{\phi}u\right\Vert _{\phi
}+\left\Vert u\right\Vert _{\phi}.$ More generally if $r>0$ is an integer,
$\mathcal{D}\left(  \overline{D_{\phi}^{r}}\right)  $ is the completion for
$\left\Vert u\right\Vert _{D_{\phi}^{r},\phi}=\left\Vert D_{\phi}%
^{r}u\right\Vert _{\phi}+\left\Vert u\right\Vert _{\phi}.$

Multiplication by $\phi$ gives a unitary $\phi:$ $\Omega_{\phi,c}^{\ast
}\rightarrow\Omega_{c}^{\ast},$ and $\phi$ induces a differential operator on
$\Omega_{c}^{\ast}$ by $\tilde{D}=\phi D_{\phi}\phi^{-1}.$ The tilde will be
used generically for operators on ordinary forms produced in this way from
operators on weighted forms. Suppose that $M$ has one component. Since
$\tilde{D}=d+\delta+$ $($zeroth order operator$)$, $\tilde{D}^{r}$ is
essentially selfadjoint by \cite{Che}. If $M$ has more than one component,
$\tilde{D}^{r}=\bigoplus_{p}\tilde{D}_{M_{p}}^{r}$ acting on $\bigoplus
_{p}\Omega_{c}^{\ast}\left(  M_{p}\right)  =$ $\Omega_{c}^{\ast}.$ By
\cite[Ex. 5.43]{Wei}, $\tilde{D}^{r}$ is essentially selfadjoint since all the
$\tilde{D}_{M_{p}}^{r}$ are. Evidently $\tilde{D}^{r}=\phi D_{\phi}^{r}%
\phi^{-1}.$ It follows that $\phi$ induces a unitary equivalence
$\overline{D_{\phi}^{r}}\rightarrow\overline{\tilde{D}^{r}}.$ Therefore
$D_{\phi}^{r}$ is also essentially selfadjoint. In particular $\left(  \bar
{D}_{\phi}\right)  ^{r}=\overline{D_{\phi}^{r}}.$

The (unreduced) $\phi$-cohomology of $M$ is defined as follows. The closures
are taken with respect to $\left\Vert \cdot\right\Vert _{\phi}.$ Let
\[
Z_{\phi}^{k}=\bigskip\newline\left\{  u\in\mathcal{D}\left(  \overline{d^{k}%
}\right)  |\overline{d^{k}}u=0\right\}  ,\;B_{\phi}^{k}=\left\{
\overline{d^{k-1}}v|v\in\mathcal{D}\left(  \overline{d^{k-1}}\right)
\right\}  .
\]
Then $H_{\phi}^{k}\left(  M\right)  =Z_{\phi}^{k}/B_{\phi}^{k}.$ Each
cohomology class has representatives which are $C^{\infty}$. Let
$\mathcal{D}^{\infty}=$ $\bigcap_{r}\mathcal{D}\left(  \bar{D}_{\phi}%
^{r}\right)  ,$ which consists of smooth forms. The inclusion $\left(
\mathcal{D}^{\infty},\bar{d}\right)  \rightarrow\left(  \mathcal{D}\left(
\bar{d}\right)  ,\bar{d}\right)  $ induces an isomorphism of cohomology
\cite[Th. 2.12]{BrLe}. Thus there is a homomorphism $S^{k}:H_{\phi}^{k}\left(
M\right)  \rightarrow H_{dR}^{k}\left(  M\right)  ,$ where the latter group is
the de Rham cohomology of $M$ based on smooth forms. For now, closures will be
understood and the bar will be suppressed. Below $H_{k}\left(  M\right)  $
denotes homology with real coefficients.

\begin{theorem}
Suppose that $H_{k}\left(  M\right)  $ is infinite dimensional. Then there is
no $\phi$ such that $S^{k}$ is surjective.
\end{theorem}

In particular, cohomology is never represented by $\phi$-harmonic forms. The
proof will be by reductio ad absurdum\textit{. }Thus assume that $S^{k}$ is surjective.

We claim that we may assume that $2k<n_{p}$ for all $p.$ This is accomplished
by taking the product of $M$ with a suitably weighted Euclidean space. The
details will be given at the end of the paper.

The homology of $M$ is certainly countably generated. Let $\left\{  \gamma
_{i}\right\}  ,i\in\mathbb{N},$ be a basis for $H_{k}\left(  M;\mathbb{Q}%
\right)  \subset H_{k}\left(  M\right)  $ which restricts to a basis of each
component. By \cite[Th. III.4]{Tho}, there are closed oriented connected
$k$-dimensional manifolds $N_{i}$, maps $g_{i}:N_{i}\rightarrow M,$ and
positive integers $m_{i}$ such that $m_{i}\gamma_{i}$ is the image of the
fundamental class of $N_{i}.$ (The statement of the cited theorem assumes that
the space is a finite polyhedron. We may triangulate $M$ and use the fact that
its homology is the direct limit of the homologies of its finite
subcomplexes.) We redefine $\gamma_{i}$ to be $m_{i}\gamma_{i},$ so that the
image is $\gamma_{i}.$

\begin{lemma}
$%
%TCIMACRO{\dcoprod \limits_{i}}%
%BeginExpansion
{\displaystyle\coprod\limits_{i}}
%EndExpansion
g_{i}$ is homotopic to an injective smooth map $%
%TCIMACRO{\dcoprod \limits_{i}}%
%BeginExpansion
{\displaystyle\coprod\limits_{i}}
%EndExpansion
f_{i}$ such that the $f_{i}\left(  N_{i}\right)  $ have disjoint closed
tubular neighborhoods $V_{i}.$
\end{lemma}

\begin{proof}
Since $2k<n_{p},$ for all $p,$ by Whitney's Embedding Theorem \cite{Whi},
$g_{1}$ is homotopic to an embedding. Let $V_{1}$ be any closed tubular
neighborhood. Assume that $f_{i}$ and $V_{i}$ have been constructed for
$i<\ell$. By transversality we can make $g_{\ell}\left(  N_{\ell}\right)  $
disjoint from $f_{i}\left(  N_{i}\right)  ,$ and then push it off $V_{i},$ for
$i<\ell.$ By the cited theorem, $g_{\ell}$ is homotopic to an embedding
$f_{\ell}$ in $M-\bigcup_{i=1}^{\ell-1}V_{i}.$ Let $V_{\ell}$ be any closed
tubular neighborhood of $f_{\ell}\left(  N_{\ell}\right)  $ in this manifold.
\end{proof}

Let the space of restrictions of elements of $\mathcal{D}\left(  D_{\phi}%
^{r}\right)  $ to $V_{i}$ be $\mathcal{D}\left(  D_{\phi}^{r}\right)  _{V_{i}%
}.$ It has the norm $\left\Vert u\right\Vert _{D_{\phi}^{r},\phi,i}$ which is
the same as $\left\Vert u\right\Vert _{D_{\phi}^{r},\phi}$ except that the
integral is evaluated on $V_{i}$. The restrictions $q_{i}:\mathcal{D}\left(
D_{\phi}^{r}\right)  \rightarrow\mathcal{D}\left(  D_{\phi}^{r}\right)
_{V_{i}}$ clearly have norm $1.$ Since the $V_{i}$ are disjoint, $\left\Vert
u\right\Vert _{D_{\phi}^{r},\phi}^{2}\geq\sum_{i}\left\Vert q_{i}u\right\Vert
_{D_{\phi}^{r},\phi,i}^{2}.$ Therefore there is a bounded operator
$r:\mathcal{D}\left(  D_{\phi}^{r}\right)  \rightarrow\widehat{\bigoplus}%
_{i}\mathcal{D}\left(  D_{\phi}^{r}\right)  _{V_{i}}$ (the Hilbert sum),
$u\rightarrow\sum_{i}q_{i}u$ \cite[Ex. 5.43]{Wei}.

We will use some properties of Sobolev spaces of $k$-forms. The basic objects
are the local spaces $W_{r,loc}^{k}\left(  M\right)  .$ For a compact
codimension zero submanifold with boundary $W$ of $M,$ $W_{r}^{k}\left(
W\right)  $ is the space of restrictions of elements of $W_{r,loc}^{k}\left(
M\right)  $ to $W$. Since $W_{0}^{k}\left(  W\right)  $ is just the $L^{2}$
space of forms, the norm will be written $\left\Vert u\right\Vert _{W}.$ See
\cite{ChPi} or \cite{Hor} for background information.

By local elliptic regularity $\mathcal{D}\left(  \tilde{D}^{r}\right)  \subset
W_{r,loc}^{\ast}\left(  M\right)  ,$ so $\mathcal{D}\left(  D_{\phi}%
^{r}\right)  \subset\phi^{-1}W_{r,loc}^{\ast}\left(  M\right)  .$ Essentially
by definition, $\phi^{-1}W_{r,loc}^{\ast}\left(  M\right)  =W_{r,loc}^{\ast
}\left(  M\right)  ,$ so $\mathcal{D}\left(  D_{\phi}^{r}\right)  \subset
W_{r,loc}^{\ast}\left(  M\right)  .$ Restricting to $V_{i},$ if $u\in
\mathcal{D}\left(  D_{\phi}^{r}\right)  ,$ then $q_{i}u\in W_{r}^{\ast}\left(
V_{i}\right)  .$ Since $\phi$ is bounded away from zero on $V_{i},$ the
unweighted and weighted graph norms satisfy
\begin{equation}
\left\Vert u\right\Vert _{D_{\phi}^{r},i}\leq L_{i}\left\Vert u\right\Vert
_{D_{\phi}^{r},\phi,i}%
\end{equation}
on $\mathcal{D}\left(  D_{\phi}^{r}\right)  _{V_{i}}$ for some constants
$L_{i}.$ Now let $U_{i}$ be a closed tubular neighborhood of $N_{i}$ contained
in the interior of $V_{i}.$ The following elliptic estimate may be found in
\cite[Th. 5.11.1]{Tay}. Let $T$ be an elliptic operator of order $r$ on $M.$
Then there is a constant $K_{i}$ such that for all $u\in W_{r}^{\ast}\left(
V_{i}\right)  ,$%
\[
\left\Vert u\right\Vert _{W_{r}^{\ast}\left(  U_{i}\right)  }\leq K_{i}\left(
\left\Vert Tu\right\Vert _{V_{i}}+\left\Vert u\right\Vert _{V_{i}}\right)  .
\]

\noindent Taking $T=D_{\phi}^{r},$ we interpret this as saying that
restriction induces a bounded operator from $\mathcal{D}\left(  D_{\phi}%
^{r}\right)  _{V_{i}}$ with the norm $\left\Vert u\right\Vert _{D_{\phi}%
^{r},i}$ to $W_{r}^{\ast}\left(  U_{i}\right)  .$ Combining this with (1),
restriction from $\mathcal{D}\left(  D_{\phi}^{r}\right)  _{V_{i}}$ with the
norm $\left\Vert u\right\Vert _{D_{\phi}^{r},\phi,i}$ to $W_{r}^{\ast}\left(
U_{i}\right)  $ has the bound $L_{i}K_{i}.$ Of course, these bounds depend on
$i.$

The final step is to evaluate forms on the fundamental classes of the $N_{i}$.
Choose $r>\tfrac{n}{2}$ so that by Sobolev's Theorem $W_{r}^{k}\left(
U_{i}\right)  $ is continuously embedded in the $C^{0}$ forms. Define a linear
functional $c_{i}$ on $W_{r}^{k}\left(  U_{i}\right)  $ by $c_{i}\left(
u\right)  =\int_{N_{i}}u.$%
\[
\left\vert c_{i}\left(  u\right)  \right\vert =\left\vert \int_{N_{i}%
}u\right\vert \leq\underset{x\in N_{i}}{Sup}\left\vert u\left(  x\right)
\right\vert Vol\left(  N_{i}\right)  \leq C_{i}\left\Vert u\right\Vert
_{W_{r}^{k}\left(  U_{i}\right)  }%
\]
for some constants $C_{i},$ so that $c_{i}$ is bounded. Then the composition
from $u\in\mathcal{D}\left(  D_{\phi}^{r}\right)  _{V_{i}}$ to $c_{i}$ has
norm less than or equal to $L_{i}K_{i}C_{i}.$ Let $\tau_{i}=\left(  L_{i}%
K_{i}C_{i}\right)  ^{-1}$and let $\mathbb{C}_{\tau_{i}}$ be $\mathbb{C}$ with
the norm $\left\Vert z\right\Vert _{\tau_{i}}=\tau_{i}\left\vert z\right\vert
$. Then the map $\widehat{\bigoplus}_{i}$ $\mathcal{D}\left(  D_{\phi}%
^{r}\right)  _{V_{i}}\rightarrow$ $\widehat{\bigoplus}_{i}\mathbb{C}_{\tau
_{i}},$ $\left(  u_{i}\right)  \rightarrow\left(  c_{i}\left(  u_{i}\right)
\right)  $ has norm $\leq1.$ Composing with $r$ gives a bounded map
$\mathcal{D}\left(  D_{\phi}^{r}\right)  \rightarrow\widehat{\bigoplus}%
_{i}\mathbb{C}_{\tau_{i}}.$

The restriction $H_{dR}^{k}\left(  M\right)  \rightarrow H_{dR}^{k}\left(
\bigcup_{i}U_{i}\right)  $ is surjective, being the (complexified) transpose
of the injection $H_{k}\left(  \bigcup_{i}U_{i}\right)  \rightarrow
H_{k}\left(  M\right)  .$ Thus a set of forms in $\mathcal{D}^{\infty}$
representing all of $H_{dR}^{k}\left(  M\right)  $ would restrict to a set
representing all of $H_{dR}^{k}\left(  \bigcup_{i}U_{i}\right)  .$ Evaluation
of elements of this group on the $N_{i}$ is well-defined, by Stokes' Theorem,
and gives an isomorphism with $\prod_{i}\mathbb{C}_{i}$ Therefore
$\widehat{\bigoplus}_{i}\mathbb{C}_{\tau_{i}}$ would contain \textit{all}
sequences of complex numbers, which is impossible. In fact, it does not
contain $\left(  \tau_{i}^{-1}\right)  .$ This completes the proof of Theorem
1 under the assumption $2k<n_{p}$ for all $p.$

We now justify this assumption. Bars will again denote closures. We will form
the product of $M$ with weight $\phi$ and $\mathbb{R}^{2k+1}$ with its usual
metric and a particular weight $\psi.$This choice satisfies the dimensional
requirement. Consider the general situation: manifolds $M_{1},$ $M_{2}$ of
dimensions $n_{1}$ and $n_{2}$ with weights $\phi$ and $\psi.$ Equip
$M_{1}\times M_{2}$ with the weight $\phi\otimes\psi.$ Let $d$ be the exterior
derivative of $M_{1}\times M_{2}$ acting on $\Omega_{\phi\otimes\psi,c}%
^{k}\left(  M_{1}\times M_{2}\right)  .$ The spaces $\mathcal{E}^{k}%
=\bigoplus_{i+j=k}\Omega_{\phi,c}^{i}\left(  M_{1}\right)  \otimes\Omega
_{\psi,c}^{j}\left(  M_{2}\right)  $ embed isometrically into $\Omega
_{\phi\otimes\psi,c}^{k}\left(  M_{1}\times M_{2}\right)  $ using the
isomorphism of exterior algebras $\Lambda_{n_{1}}\hat{\otimes}\Lambda_{n_{2}%
}\cong\Lambda_{n_{1}+n_{2}},$ $\hat{\otimes}$ the graded tensor product.
Denote the exterior derivatives of $M_{1},$ $M_{2},$ and $M_{1}\times M_{2}$
by $d_{1},$ $d_{2},$ and $d.$ The restriction of $d$ to $\mathcal{E}^{k}$ is
given by $\bigoplus_{i+j=k}\left(  d_{1}^{i}\otimes I+\left(  -1\right)
^{i}I\otimes d_{2}^{j}\right)  .$ There is a homomorphism $H_{\phi}^{i}\left(
M_{1}\right)  \otimes H_{\psi}^{j}\left(  M_{2}\right)  \rightarrow
H_{\phi\otimes\psi}^{k}\left(  M_{1}\times M_{2}\right)  .$ Choose
representatives $u,$ $v$ of given cohomology classes which are in
$\mathcal{D}^{\infty}$. We then check directly that $u\otimes v$ is a smooth
closed form in $\mathcal{D}\left(  \bar{d}\right)  .$ Its class in
$H_{\phi\otimes\psi}^{k}\left(  M_{1}\times M_{2}\right)  $ is independent of
the choices. The following diagram is obviously commutative.%
\[%
\begin{tabular}
[c]{ccc}%
$\bigoplus_{i+j=k}\left(  H_{\phi}^{i}\left(  M_{1}\right)  \otimes H_{\psi
}^{j}\left(  M_{2}\right)  \right)  $ & $\overset{\bigoplus\left(  S_{1}%
^{i}\otimes S_{2}^{j}\right)  }{\longrightarrow}$ & $\bigoplus_{i+j=k}\left(
H_{dR}^{i}\left(  M_{1}\right)  \otimes H_{dR}^{j}\left(  M_{2}\right)
\right)  $\\
$\downarrow$ &  & $\downarrow$\\
$H_{\phi\otimes\psi}^{k}\left(  M_{1}\times M_{2}\right)  $ & $\overset
{S}{\longrightarrow}$ & $H_{dR}^{k}\left(  M_{1}\times M_{2}\right)  .$%
\end{tabular}
\ \ \ \ \ \ \ \ \
\]
On the right we use any smooth representatives of the classes. From \cite[Prop
II.9.12]{BoTu} we have that the right arrow is an isomorphism provided that
$H_{dR}^{\ast}\left(  M_{2}\right)  $ is finitely generated. (These authors
use an equivalent description of the map as $\pi_{1}^{\ast}u\wedge\pi
_{2}^{\ast}v.)$ (The left arrow is also an isomorphism if $H_{\psi}^{\ast
}\left(  M_{2}\right)  $ is finitely generated \cite[Th. 2.14]{BrLe}, but this
isn't needed.)

Let $M_{1}=M$ and $M_{2}=\mathbb{R}^{2k+1}.$ Let $\psi=e^{-\left\vert
x\right\vert },$ smoothed near the origin. An application of \cite[Th. 3.1,
Rem. p.161]{Yeg} shows that $s_{2}$ is an isomorphism. The diagram reduces to
\[%
\begin{tabular}
[c]{ccc}%
$H_{\phi}^{k}\left(  M\right)  $ & $\overset{s_{1}^{k}}{\longrightarrow}$ &
$H_{dR}^{k}\left(  M\right)  $\\
$\downarrow$ &  & $\downarrow\cong$\\
$H_{\phi\otimes\psi}^{k}\left(  M\times\mathbb{R}^{2k+1}\right)  $ &
$\overset{s^{k}}{\longrightarrow}$ & $H_{dR}^{k}\left(  M\times\mathbb{R}%
^{2k+1}\right)  .$%
\end{tabular}
\ \ \ \ \ \
\]
If the top arrow were surjective, the bottom one would be as well,
contradicting the statement already proved.


\begin{thebibliography}{99}                                                                                               %


\bibitem {Bor}Borel, A. Stable and $L_{2}$-cohomology of arithmetic groups.
\textit{Bull. Amer. Math. Soc. (N.S.)} \textbf{3} (1980), no. 3, 1025--1027.

\bibitem {BoTu}Bott, R.; Tu, L. W. Differential forms in algebraic topology.
Graduate Texts in Mathematics, 82. \textit{Springer-Verlag, New York-Berlin}, 1982.

\bibitem {BrLe}Br\"{u}ning, J.; Lesch, M. Hilbert complexes. \textit{J. Funct.
Anal. }\textbf{108} (1992), no. 1, 88--132.

\bibitem {Bul}Bullock, S. S. Gaussian weighted unreduced $L_{2}$ cohomology of
locally symmetric spaces. \textit{New York J. Math}. \textbf{8} (2002),
241--256 (electronic).

\bibitem {Bue}Bueler, E. L. The heat kernel weighted Hodge Laplacian on
noncompact manifolds. \textit{Trans. Amer. Math. Soc}. \textbf{351} (1999),
no. 2, 683--713.

\bibitem {Car}Carron, G. A Counter Example to the Bueler's
Conjecture\textit{.} \textit{www.arxiv.org/abs/math/0509550.}

\bibitem {ChPi}Chazarain, J.; Piriou, A. Introduction to the theory of linear
partial differential equations. Translated from the French. Studies in
Mathematics and its Applications, 14. \textit{North-Holland Publishing Co.,
Amsterdam-New York}, 1982.

\bibitem {Che}Chernoff, P. R. Essential self-adjointness of powers of
generators of hyperbolic equations. \textit{J. Functional Analysis}
\textbf{12} (1973), 401--414.

\bibitem {Dod}Dodziuk, J. de Rham-Hodge theory for $L_{2}$-cohomology of
infinite coverings. \textit{Topology} \textbf{16} (1977), no. 2, 157--165.

\bibitem {Hor}H\"{o}rmander, L. Linear partial differential operators.
\textit{Springer Verlag, Berlin-New York}, 1976.

\bibitem {Mil}Miller, J. G. The Euler characteristic and finiteness
obstruction of manifolds with periodic ends.\textit{ Asian J. Math.
}\textbf{10 }(2006), no. 4, 679-714.

\bibitem {Tay}Taylor, M. E. Partial differential equation\textit{s. I. }Basic
theory. Applied Mathematical Sciences, 115. \textit{Springer-Verlag, New
York}, 1996.

\bibitem {Tho}Thom, R. Quelques propri\'{e}t\'{e}s globales des
vari\'{e}t\'{e}s diff\'{e}rentiables. \textit{Comment. Math. Helv.}
\textbf{28}, (1954). 17--86.

\bibitem {Wei}Weidmann, J. Linear operators in Hilbert spaces. Translated from
the German by Joseph Sz\"{u}cs. Graduate Texts in Mathematics, 68.
\textit{Springer-Verlag, New York-Berlin}, 1980.

\bibitem {Whi}Whitney, H. Differentiable manifolds. \textit{Ann. of Math. (2})
\textbf{37} (1936), no. 3, 645--680.

\bibitem {Yeg}Yeganefar, N. Sur la $L_{2}$-cohomologie des vari\'{e}t\'{e}s
\`{a} courbure n\'{e}gative. \textit{Duke Math. J.} \textbf{122} (2004), no.
1, 145--180.

\bibitem {Zuc}Zucker, S. $L_{2}$ cohomology of warped products and arithmetic
groups. \textit{Invent. Math.} \textbf{70} (1982/83), no. 2, 169--218.
\end{thebibliography}
\end{document}